\documentclass[12pt]{article}
\usepackage{amssymb}
\usepackage{amsmath}
\usepackage{amstext}
\usepackage{amsfonts}
\usepackage{amsthm}
\usepackage{amscd}
\numberwithin{equation}{section}

\swapnumbers
\theoremstyle{plain}
\newtheorem{thm}{Theorem}[section]
\newtheorem{lem}[thm]{Lemma}
\newtheorem{prop}[thm]{Proposition}
\newtheorem{cor}[thm]{Corollary}

\theoremstyle{remark}
\newtheorem{rem}[thm]{Remark}

\def\cO{\mathcal{O}}

\def\map#1{\ \smash{\mathop{\longrightarrow}\limits^{#1}}\ }

\def\LL{\mathcal{L}}
\def\XX{\mathcal{X}}
\def\S{\mathfrak{S}}
\def\MM{\mathrm{M}}
\def\NN{\mathbb{N}}
\def\GG{\mathbb{G}}

\def\ZZ{\mathbb{Z}}

\def\QQ{\mathbb{Q}}
\def\FF{\mathbb{F}}
\def\AA{\mathbb{A}}

\def\HH{\mathcal{H}}
\def\VV{\mathcal{V}}
\def\WW{\mathcal{W}}
\def\II{\mathcal{I}}

\def\div{\mathrm{div}}
\def\deg{\mathrm{deg} \:}

\def\ker{\mathrm{ker} \:}
\def\pic{\mathrm{Pic}}

\def\kum{\mathrm{Kum} \:}
\def\dim{\mathrm{dim} \:}
\def\det{\mathrm{det} \:}

\def\Spec{\mathrm{Spec}}

\def\lra{\longrightarrow}

\def\ra{\rightarrow}

\def\lms{\longmapsto}

\def\kum{\mathrm{Kum}}
\def\GL{\mathrm{GL}}

\def\pp{\mathbb{P}}

\def\M{\mathcal{M}}

\def\EE{\mathcal{E}}
\def\JJ{\mathcal{J}}

\title{The Frobenius map, rank $2$ vector bundles and Kummer's quartic 
surface in characteristic $2$ and $3$}

\author{Yves Laszlo and Christian Pauly}
\begin{document}
\maketitle


\section{Introduction}
Let $X$ be a smooth projective curve of genus $2$ defined over an algebraically
closed field $k$ of characteristic $p>0$. The moduli space $\MM_X$ of
semi-stable rank $2$ vector bundles with fixed trivial
determinant is isomorphic to the linear system $|2\Theta| \cong \pp^3$
over $\pic^1(X)$ and the $k$-linear relative Frobenius map $F: X \ra X_1$
induces by pull-back a rational map (the Verschiebung)
\begin{equation} \label{commdiagmod}
\begin{CD}
\MM_{X_1} @>V>> \MM_X \\
@VDVV  @VVDV \\
|2\Theta_1| @>\tilde{V}>> |2\Theta|
\end{CD}
\end{equation}
The vertical maps $D$ are isomorphisms and the Verschiebung $V: E \mapsto
F^* E$ coincides via $D$ with a rational map $\tilde{V}$ given by 
polynomial equations of degree $p$ (Proposition \ref{degreeV}).
The Kummer surfaces $\kum_X$ and $\kum_{X_1}$ are 
canonically contained in the
linear systems $|2\Theta|$ and $|2\Theta_1|$
and coincide with the semi-stable boundary of the moduli spaces
$\MM_X$ and $\MM_{X_1}$. Moreover $\tilde{V}$
maps $\kum_{X_1}$ onto  $\kum_X$.

\bigskip

Our interest in the diagram \eqref{commdiagmod} comes from
questions related to the action of the Frobenius map on vector
bundles like e.g. surjectivity of $V$, density of Frobenius-stable
bundles, loci of Frobenius-destabilized bundles (see \cite{lp}). These
questions are largely open when the rank of the bundles, the genus
of the curve or the characteristic of the field are arbitrary. 
In \cite{lp} we made use of the exceptional isomorphism 
$D: \MM_X \ra |2\Theta|$ in the genus $2$, rank $2$ case 
and determined the equations of $\tilde{V}$
when $X$ is an ordinary curve and $p=2$, which allowed us to answer the
above mentioned questions. In this
paper we obtain the equations of $\tilde{V}$ in two more cases:
\begin{itemize}
\item[(1)] $p = 2$ and $X$ non-ordinary with Hasse-Witt
invariant equal to 1,
\item[(2)] $p = 3$ and any $X$.
\end{itemize}

In case (1) we consider a family $\XX$ of genus $2$ curves 
parametrized by a discrete valuation ring with ordinary generic fibre
$\XX_\eta$ and special fibre isomorphic to $X$. We obtain the
equations of $\tilde{V}$ for $X$ (Theorem \ref{mainthm})
by specializing the quadrics $P_{ij}$ defining the Verschiebung
$V_\eta : \MM_{\XX_{1\eta}} \lra \MM_{\XX_\eta}$
associated to $\XX_\eta$ (section 5). In order to determine
the limit of the $P_{ij}$'s we use the explicit formulae 
(Proposition \ref{coeffkum}) of the coefficients of the $P_{ij}$'s,
which coincide with the coefficients of Kummer's quartic
surface $\kum_{\XX_\eta}$, in terms of the coefficients of an affine equation
of the ordinary curve $\XX_\eta$.
As in the ordinary case we easily deduce from the equations 
of $\tilde{V}$  a full description of the
Verschiebung $V$ (Proposition \ref{descr}).

\bigskip

In case (2) we show that the cubic equations of $\tilde{V}$ are 
given by the polar equations of a Kummer surface $S \subset |2\Theta_1|$
(Theorem \ref{theochar3}). Moreover $S$ is isomorphic to
$\kum_X$ and the $16$ nodes of $S$ correspond to the $16$ base points
of $\tilde{V}$. We deduce that $V$ is surjective and of degree $11$
(Corollary \ref{corchar3}).

\bigskip

In the appendix we show that over any smooth curve $X$ of genus $g \geq 2$ 
defined over an algebraically closed field of characteristic $p>0$ 
and for any integer $r \geq 2$, there exist  
Frobenius-destabilized bundles of rank $r$, i.e., 
semi-stable bundles $E$ such that $F^*E$ is not semi-stable
(Theorem \ref{exisbp}).

\bigskip

We thank M. Raynaud for helpful discussions.

\section{Preliminaries on genus $2$ curves in characteristic $2$}

We consider a smooth {\em ordinary} curve $X$ of genus $2$ defined 
over a field $K$ of characteristic $2$. We denote by $\overline{K}$ its
algebraic closure. Let $k(X)$ be the
function field of $X$ and $\omega_X$ the canonical bundle over $X$.

\subsection{Weierstrass points}
After taking a finite extension of $K$  and applying an automorphism
of $\pp^1_K$ we can assume that the three 
Weierstrass points of $X$ are $0,1$ and $\infty$. 
We consider the birational Abel-Jacobi map 
\begin{equation} \label{abel}
AJ: \ S^2 X \lra JX, \qquad P_1 + P_2 \lms \cO_X(P_1 +P_2)
\otimes \omega_X^{-1}. 
\end{equation}
We observe that the three nonzero elements
$[0],[1],[\infty]$ of $JX[2]$ are $K$-rational and are given by 
\begin{equation}\label{twotorsion}
[0]:= AJ(1+\infty), \qquad
[1]:= AJ(0+\infty), \qquad
[\infty]:= AJ(0+1).
\end{equation}
For later use we mention that the sheaf of locally exact differential
forms (see \cite{ray} section 4) equals
$$ B = \cO_X(0+1+\infty)\otimes \omega_X^{-1}. $$
We recall that $B$ is a theta-characteristic of $X$, i.e.,
$B^{\otimes 2} \cong \omega_X$.

\subsection{Level $2$ structure}

A level $2$ structure is an isomorphism 
$\psi: JX[2] \map{\sim} \FF_2^2$. Note that two level
$2$ structures differ by an automorphism of $\FF^2_2$, i.e.,
an element of $\GL(2,\FF_2) = \S_3$, where $\S_3$ is the symmetric group.
It is well-known that a level $2$ structure $\psi$ is
equivalent to an ordering of the Weierstrass points of $X$.
We refer e.g. to \cite{do} page 141 for the characteristic zero case,
which can easily be adapted to the characteristic two case.
Because of the choices involved in the degeneration
of $X$ to a non-ordinary curve (section 5.1), we consider 
the ordering $1,\infty,0$ of the Weierstrass points.
With the notation of \eqref{twotorsion}  the corresponding level
$2$ structure $\psi$ is given by
\begin{equation} \label{leveltwopsi} 
\psi([0]) = (1,0), \qquad \psi([1]) = (0,1), \qquad \psi([\infty]) = (1,1).
\end{equation}

A level $2$ structure $\psi$ allows us to construct (\cite{lp} section 2)
the theta basis $\{ X_g \}_{g \in \FF^2_2}$ of the space $H^0(JX,2\Theta)$.
We denote the four sections by $X_B,X_{0},X_{1},X_{\infty}$ and
introduce the rational functions $Z_{\bullet} \in k(JX)$ defined by
\begin{equation} \label{rattheta}
Z_{0} = \frac{X_{0}}{X_B}, \qquad
Z_{1} = \frac{X_{1}}{X_B}, \qquad
Z_{\infty} = \frac{X_{\infty}}{X_B}.
\end{equation}
We recall that $X_B$ is the theta function with associated zero
divisor $2\Theta_B$, where
$$ \Theta_B := \{ L \in JX \ | \ h^0(X,L \otimes B) \geq 1 \}. $$
 
\subsection{Birational models}

Let $X$ be an ordinary smooth curve of genus $2$ defined over
$\overline{K}$ and $\psi$ a level $2$ sturcture.
It follows from \cite{lange} page 28  and \cite{bhosle}
Proposition 1.4 that the pair $(X,\psi)$ is {\em uniquely}
represented by an affine equation of the form
\begin{equation} \label{affmodbhosle}
y^2 + x(x+1) y = x(x+1)(ax^3 + (a+b)x^2 +cx + c), 
\end{equation} 
with $a,b,c \in \overline{K}$. Moreover if $X$ is defined over $K$,
the coefficients $a,b,c$ lie in a finite extension of $K$.
The next lemma is an immediate consequence of \cite{bhosle} 
Proposition 1.5.

\begin{lem}
The curve $X$ defined by the equation \eqref{affmodbhosle} is
smooth if and only if $abc \not= 0$.
\end{lem}

Let $\widetilde{\M}_3$ denote the moduli space parametrizing pairs $(X,\psi)$
of smooth  ordinary genus $2$ curves $X$ defined over $\overline{K}$ 
equipped with a  level $2$ structure. It follows from the previous
remark that $\widetilde{\M}_3$ is an affine variety,
$$ \widetilde{\M}_3 = \Spec \ \overline{K}[a,b,c,\frac{1}{abc}]. $$

\noindent
Fixing the curve $X$ the symmetric group $\S_3$ acts naturally 
on the level $2$ structures $\psi$. It can be shown that this $\S_3$-action
on $\widetilde{\M}_3$ coincides with the permutation action of 
$\S_3$ on the coefficients $a,b,c$.

\subsection{Normal form} \label{nf}

We introduce the rational function $Y \in k(X)$ 
defined by $Y = \frac{y}{x(x+1)}$.
Then \eqref{affmodbhosle} becomes 
\begin{equation} \label{normalform}
Y^2 + Y = R(x), \qquad \text{with} \ R(x) = 
\frac{ax^3 + (a+b)x^2 + cx +c}{x(x+1)}.
\end{equation}
We also observe that, given a polynomial $S(x)$, the
involution $i_S: \AA^2_k \ra \AA^2_k$, $(x,y) \mapsto (x,y+x(x+1)S(x))$
transforms the equation \eqref{normalform} of $X$ into $Y^2 + Y 
= R +S^2 + S$. 

\subsection{Kummer's quartic equation}

For a pair $(X,\psi)$ it has been shown in \cite{lp} Proposition 4.1 that
there exist constants $\lambda_{0},\lambda_{1},\lambda_{\infty}
\in \overline{K}$ such that the following equality holds in $k(JX)$,
\begin{equation}\label{kumeq}
\lambda_{0}^2 (Z_{0}^2 + Z_{1}^2Z_{\infty}^2) +
\lambda_{1}^2 (Z_{1}^2 + Z_{0}^2Z_{\infty}^2) +
\lambda_{\infty}^2 (Z_{\infty}^2 + Z_{0}^2Z_{1}^2) +
\lambda_{0}\lambda_{1}\lambda_{\infty}Z_{0}Z_{1}
Z_{\infty} = 0.
\end{equation}
The constants $\lambda_0,\lambda_1,\lambda_\infty$ are related via $\psi$
\eqref{leveltwopsi} to the $\{ \lambda_g \}_{g \in \FF_2^2}$ used in 
\cite{lp} Proposition 4.1 as follows: $\lambda_{0} = 
\frac{\lambda_{10}}{\lambda_{00}}$, $\lambda_{1} = 
\frac{\lambda_{01}}{\lambda_{00}}$, $\lambda_{\infty} = 
\frac{\lambda_{11}}{\lambda_{00}}$.

\section{The coefficients $\lambda_\bullet$ of Kummer's quartic surface 
$\kum_X$ for ordinary $X$}

\begin{prop} \label{coeffkum}
Given a curve $X$ with a level $2$ structure $\psi$ represented by an 
affine equation \eqref{affmodbhosle} with $a,b,c \in K$. Then the 
coefficients of the equation \eqref{kumeq} of its Kummer surface 
$\kum_X$ are
$$ \lambda^2_{0} = \frac{1}{ab}, \qquad
 \lambda^2_{1} = \frac{1}{ac}, \qquad
 \lambda^2_{\infty} = \frac{1}{bc}. $$
Let $\{ x_g \}$ be the dual basis of the theta basis $\{ X_g \}_{g \in
\FF^2_2}$ of $|2\Theta|= \pp^3$. Then the homogeneous equation of $\kum_X$ is
$$ c (x^2_{00}x^2_{10} + x^2_{01}x^2_{11}) +
b (x^2_{00}x^2_{01} + x^2_{10}x^2_{11})  
+ a (x^2_{00}x^2_{11} + x^2_{01}x^2_{10})  +
x_{00}x_{01}x_{10}x_{11} = 0. $$
\end{prop}

The idea of the proof is to consider the pull-back of the
rational functions $Z_{\bullet}$ \eqref{rattheta} by the
Abel-Jacobi map \eqref{abel} to the symmetric product
$S^2 X$ and to do the computations in the function field $k(S^2 X)
\hookrightarrow k(X \times X)$. Since $X$ is given by the equation 
\eqref{affmodbhosle}, we have
natural coordinates $x_1,y_1$ and $x_2,y_2$ on $X \times X$.
For notational convenience, we also use $Y_i = \frac{y_i}{x_i(x_i+1)}$
for $i=1,2$.

\noindent
We use two lemmas.

\begin{lem} \label{rellem}
Suppose that there exist polynomials $A,B \in K[x_1,x_2]$, which
satisfy
\begin{equation} \label{rlem}
(Y_1 + Y_2) A(x_1,x_2) = B(x_1,x_2).
\end{equation}
Then $A=B=0$.
\end{lem}

\begin{proof}
Squaring the relation \eqref{rlem} and using \eqref{normalform} leads
to the equation
$$(Y_1 + Y_2) A^2 + (R(x_1) + R(x_2))A^2 + B^2 = 0. $$
Applying again \eqref{rellem}, the first term transforms into $AB$. Clearing
denominators, we arrive at a polynomial equation, which only holds
if $A=B= 0$; e.g. take the total degree of $A$ and $B$.
\end{proof}

\begin{lem} \label{pullback}
The pull-back by the Abel-Jacobi map $AJ$ of the rational function
$Z_\infty \in k(JX)$ equals
$$ 
AJ^*(Z_\infty) = \alpha_\infty \frac{P(x_1,x_2)}{(Y_1 + Y_2)^2}, 
\qquad \text{with} \ P(x_1,x_2) = \frac{(x_1 + x_2)^2}{x_1x_2(x_1+1)(x_2 +1)} 
$$
Similary we have
$$
AJ^*(Z_0) = \alpha_0 \frac{P(x_1,x_2)}{(Y_1 + Y_2)^2}  
x_1 x_2 , \qquad
AJ^*(Z_1) = \alpha_1 \frac{P(x_1,x_2)}{(Y_1 + Y_2)^2}(x_1 +1)(x_2 +1),
$$
for some nonzero constants $\alpha_0,\alpha_1,\alpha_\infty \in K$.
\end{lem}

\begin{proof}
The first equality follows immediately from Theorem 2 \cite{ag}
and the other two from Proposition 5 \cite{ag}.
\end{proof}

\begin{proof}[Proof of Proposition \ref{coeffkum}]

We write $Q= \frac{P(x_1,x_2)}{(Y_1 + Y_2)^2}$. Using Lemma
\ref{pullback} the pull-back to $S^2 X$ of the equation \eqref{kumeq}
equals
\begin{align*}
\lambda_0^2 \left[\alpha^2_0 x_1^2 x_2^2 Q^2 + \alpha_1^2 \alpha_\infty^2 
(x_1+1)^2(x_2+1)^2 Q^4 \right] +
\lambda_1^2 \left[\alpha^2_1 (x_1+1)^2(x_2+1)^2 Q^2 + 
\alpha_0^2 \alpha_\infty^2 
x_1^2x_2^2 Q^4 \right] + \\
\lambda_\infty^2 \left[\alpha^2_\infty Q^2 + \alpha_0^2 \alpha_1^2 
x_1^2x_2^2(x_1+1)^2(x_2+1)^2 Q^4 \right]  
+ \lambda_0 \lambda_1 \lambda_\infty \left[ \alpha_0 \alpha_1 \alpha_\infty
x_1 x_2 (x_1+1) (x_2+1) Q^3 \right] = 0
\end{align*}

\noindent
We divide by $Q^2$ and multiply by $(Y_1 + Y_2)^4$,
\begin{align} \label{eq}
(Y_1 + Y_2)^4 \left[\lambda_0^2 \alpha_0^2 x_1^2 x_2^2 + 
\lambda_1^2 \alpha_1^2 (x_1 +1)^2(x_2+1)^2 + \lambda^2_\infty \alpha^2_\infty
\right] + \nonumber \\ 
(Y_1 + Y_2)^2 \left[\lambda_0 \lambda_1 \lambda_\infty \alpha_0 \alpha_1 
\alpha_\infty P x_1 x_2 (x_1+1) (x_2+1) \right] + S = 0.  
\end{align}
where $S$ is the sum of the remaining terms (not containing $Y_1,Y_2$).
After taking the square root (note that the entire expression is a square
of a polynomial in the $x_i$'s and $Y_i$'s), applying \eqref{normalform}
and Lemma \ref{rellem}, we obtain that the coefficients of $(Y_1+Y_2)^2$
and $(Y_1+Y_2)^4$ are the same, i.e.,
$$ \lambda_0^2 \alpha_0^2 x_1^2 x_2^2 + 
\lambda_1^2 \alpha_1^2 (x_1 +1)^2(x_2+1)^2 + \lambda^2_\infty \alpha^2_\infty
 = \lambda_0 \lambda_1 \lambda_\infty \alpha_0 \alpha_1 
\alpha_\infty P x_1 x_2 (x_1+1)(x_2+1). $$
An easy computation shows that this equality only holds if
$$ \lambda_0 \alpha_0 = \lambda_1 \alpha_1 = \lambda_\infty \alpha_\infty
= 1. $$
Now we replace the $\alpha$'s and the sum $S$ by their expressions in the
square root of the equation \eqref{eq}
\begin{align*}
(Y_1 + Y_2)^2 (x_1 + x_2) + (Y_1 + Y_2) (x_1 + x_2) + \\
P  \left[ \frac{\lambda_0}{\lambda_1 \lambda_\infty} (x_1+1)(x_2+1)
 + \frac{\lambda_1}{\lambda_0 \lambda_\infty} x_1 x_2 +
\frac{\lambda_\infty}{\lambda_0 \lambda_1} x_1 x_2 (x_1+1) (x_2+1)
\right] = 0.
\end{align*}
We introduce the constants $\mu_0, \mu_1, \mu_\infty \in K$ defined 
by  $\mu_\infty = \frac{\lambda_\infty}
{\lambda_0 \lambda_1}, \mu_0 = \frac{\lambda_0}{\lambda_1 \lambda_\infty},
\mu_1 = \frac{\lambda_1}{\lambda_0 \lambda_\infty}$. The previous
equality becomes after replacing $P$ by its expression and dividing
by $(x_1 +x_2)$
$$
\left[ Y_1^2 + Y_1 + \frac{\mu_0}{x_1} + \frac{\mu_1}{x_1 + 1}
+ \mu_\infty x_1 \right] +
\left[ Y_2^2 + Y_2 + \frac{\mu_0}{x_2} + \frac{\mu_1}{x_2 + 1}
+ \mu_\infty x_2 \right]  = 0. 
$$
This equation holds in $k(X \times X)$ and since the variables $(x_1,Y_1)$ and
$(x_2,Y_2)$ are separated, each of the two terms equals zero. So we can drop
the indices and we obtain an equation
\begin{equation} \label{eq2}
Y^2 + Y = \mu_\infty x + \frac{\mu_0}{x}  + 
\frac{\mu_1}{x + 1} = \frac{\mu_\infty x^3 + \mu_\infty x^2 + 
(\mu_0 + \mu_1) x + \mu_0}{x(x+1)}, 
\end{equation}
which has to be equivalent (after applying an automorphism of
$\AA_{\overline{K}}^2$) to the normal form \eqref{normalform} 
of the equation of $X$.
The automorphism is given by $i_S$ (see section \ref{nf}) 
with $S(x) = s \in \overline{K}$ satisfying $s^2 +s = \mu_1$.
Hence \eqref{eq2} is equivalent to $Y^2 + Y = R(x)$ with 
$R(x) =  \frac{\mu_\infty x^3 + (\mu_\infty + \mu_1)x^2 + \mu_0x + \mu_0}
{x(x+1)}$. Hence by uniqueness of the normal form, we obtain
$a= \mu_\infty,b=\mu_1, c= \mu_0$ and therefore also the relations claimed
in the proposition.
\end{proof}






\section{Deformation of genus $2$ curves}

Let $X/k$ be a smooth curve with Hasse-Witt invariant equal to $1$. 
We introduce a family
$\XX$ over $R = k[[t]]$ such that the special fibre $\XX_0$ is 
isomorphic to $X$ and the generic fibre $\XX_\eta$ is an ordinary
genus $2$ curve. Here $\eta$ (resp. $0$) is the generic (resp. closed)
point of $\Spec(R)$. In section 5.1 we give an example of a family 
$\XX$ with given special fibre $X$. Let $\JJ \XX$ be its associated Jacobian
scheme over $\Spec(R)$.

\subsection{$2$-divisible groups}

Let $\JJ\XX(2)$ be the $2$-divisible group of $\JJ\XX$, which is finite
and flat over $\Spec(R)$. We consider the canonical exact sequence 
\begin{equation} \label{esconet}
0 \lra \JJ\XX(2)^0 \lra \JJ\XX(2) \lra \JJ\XX(2)^{et} \lra 0,
\end{equation}
where $\JJ\XX(2)^0$ (resp. $\JJ\XX(2)^{et})$ is a connected
(resp. \'etale) $2$-divisible group. Taking again the connected
component of the Cartier dual of $\JJ\XX(2)^0$ we obtain a 
filtration
$$\JJ\XX(2)^{00} \subset \JJ\XX(2)^0 \subset \JJ\XX(2),$$
with quotients given by the $2$-divisible groups
$$ \JJ\XX(2)/\JJ\XX(2)^0 = \JJ\XX(2)^{et} \cong \QQ_2/\ZZ_2, \qquad
\JJ\XX(2)^0/ \JJ\XX(2)^{00} \cong \GG_m(2).$$
The $2$-divisible group $\JJ\XX(2)^{00}$ is self-dual,
of dimension $1$ and of height $2$. Because of the uniqueness
of $2$-divisible groups over $k$ with these properties 
(see e.g. \cite{de} Examples page 93), the special fibre $\JJ\XX(2)^{00}_0$ 
$(= \JJ\XX(2)^{00} \otimes_R k)$ is isomorphic to the $2$-divisible
 group associated
to the supersingular elliptic curve $E^{ss}/k$. We recall that there
exists a unique (up to isomorphism) supersingular curve $E^{ss}$, which
is defined by $j=0$. Therefore by a theorem of Serre-Tate \cite{ka}, 
there exists an elliptic
curve $\EE_\XX$ over $\Spec(R)$ such that $(\EE_\XX)_0 \cong E^{ss}$ and the 
associated $2$-divisible group $\EE_\XX(2)$ is isomorphic to 
$\JJ\XX(2)^{00}$ over  $\Spec(R)$.

\subsection{Deformation of elliptic curves}

In this section we compute the  linear action of the $2$-torsion
point of $\EE$ on the space of second order theta functions 
$H^0(\EE, 2\Theta)$ for a family of elliptic curves $\EE/\Spec(R)$ 
with supersingular special fibre $\EE_0 \cong E^{ss}$ and 
ordinary generic fibre
$\EE_\eta$.

\subsubsection{Addition on an ordinary elliptic curve}

Let $E$ be an ordinary elliptic curve defined over a field $K$ by the
homogeneous equation
\begin{equation} \label{eqcub}
Y^2Z + a_1 XYZ = X^3 + a_2 X^2Z + a_4 XZ^2
\end{equation}
with $a_1,a_2,a_4 \in K$ and $a_1 \not= 0$. We take as origin the inflection
point $\infty$ with projective coordinates $(0:1:0)$. The projection with
center $\infty$ gives a $2:1$ morphism $E \map{\pi} \pp^1_{K}$, 
with $X,Z$ projective
coordinates on $\pp^1_{K}$. The Abel-Jacobi map $E \ra JE$, $e \mapsto
\cO_E(e - \infty)$ identifies $E$ with $JE$. Under this identification
the theta divisor $\Theta_B$, associated to the canonical 
theta-characteristic $B= \cO_E(P - \infty) \in JE[2]$, becomes 
the $2$-torsion point $P$ with projective coordinates $(0:0:1)$.
Moreover, using this identification, we have $\cO_{JE}(2\Theta) =
\cO_E(2P) = \pi^* \cO_{\pp^1}(1)$.

The point $B \in JE[2]$ induces a linear involution, denoted by $g$, on 
the space
\begin{equation} \label{thetaell} 
W = H^0(JE, \cO(2\Theta)) =  H^0(E, \pi^* \cO_{\pp^1}(1)),
\end{equation}
such that for all nonzero $s \in W$ we have $T_B \div(s) = 
\div(g.s)$. Here $T_B$ denotes translation in $JE$ by the point $B$.
The space $W$ has two distinguished bases: first the coordinate
functions $\{X,Z\}$ and secondly the Theta basis $\{X_0,X_1\}$
(see \cite{lp} section 2).
Since the canonical section $X_0 \in W$ (associated to 
the divisor $\Theta_B$) is proportional to $X$, there exists a nonzero
$a \in K$ such that
$$ g. X = a Z, \qquad g. Z = a^{-1} X. $$
In order to determine $a$ in terms of the coefficients
$a_i \in K$, we choose one of the two points on $E$ with projective
coordinates of the form $(1:Y:1)$ and call it $A$. By construction
we have $A \in \div(X+Z)$ and, after applying $T_B$, we obtain
$A + P \in T_B \div(X+Z)$. Since 
$$ T_B \div(X+Z) = \div(g.X + g.Z) = \div(aZ + a^{-1} X),$$
we deduce that 
$$ \left( \frac{X}{Z} \right)(A+P) = a^2.$$
Now the addition formula for elliptic curves (see e.g. \cite{sil}
page 59) implies that $\left( \frac{X}{Z} \right)(A+P) = a_4$.
Hence $a = (a_4)^{1/2}$ and the Theta basis of $W$ is given by
\begin{equation} \label{heisactell}
X_0 = X, \qquad X_1 = g.X_0 = (a_4)^{1/2} Z.
\end{equation}

\subsubsection{An example}

We consider the family of elliptic curves $\EE$ over $\Spec(R)$ 
defined by the homogeneous equation
\begin{equation} \label{eqcubfam}
V^2Z + t^4 UVZ + VZ^2 = U^3 + UZ^2 
\end{equation}
with origin $\infty = (0:1:0)$. The generic fiber $\EE_\eta$ is an
ordinary elliptic curve over $\Spec(K)$, with $K = k((t))$, and the
special fibre is supersingular,i.e., $\EE_0 \cong E^{ss}$. 
The $2$-torsion point  $P_\eta$ of $\EE_\eta$ has projective coordinates
$$P_\eta = t^6 (u_0:v_0:1),$$
with $u_0 = t^{-4}$, $v_0 = t^{-2} + t^{-6}$, which specializes
to $\infty_0 \in \EE_0$. 

The $R$-module $H^0(\EE, 2\Theta)$ is free, of rank $2$, and an 
$R$-basis is given by $\{ U,Z \}$. In order to compute the linear action
$g$ of the $2$-torsion point $P \in \EE$ on $H^0(\EE,2\Theta)$, we
consider the generic fibre $\EE_\eta$. The change of variables
$X = U + u_0 Z$ and $Y = V + v_0 Z$ transforms equation
\eqref{eqcubfam} into \eqref{eqcub}, with $a_1 = t^4$, $a_2 = u_0$,
and $a_4 = 1 + u_0^2 + t^4 v_0 = 1 + t^{-8} + t^2 + t^{-2}$.
With the notation of section 4.2.1 we find 
$a = t^{-4} (1+t^3 + t^4 + t^5) = (a_4)^{1/2}$.
Therefore the action of $g$ on $H^0(\EE,2\Theta)$ is given by 
the formulae
\begin{align*}
g.U = \frac{1}{1+t^3 + t^4 + t^5} U + 
\frac{t^2 + t^4 + t^6}{1+t^3 + t^4 + t^5} Z, \\
g.Z = \frac{t^4}{1+t^3 + t^4 + t^5} U + 
\frac{1}{1+t^3 + t^4 + t^5} Z, 
\end{align*}
and, using \eqref{heisactell}, the theta functions $X_0$ and
$X_1$ can be expressed in the $R$-basis $\{ U,Z \}$ as 
follows
$$
\begin{array}{ccl}
X_0 & = & t^4 U + Z, \\
X_1 & = & g.X_0 = (1+t^3 + t^4 + t^5)Z.  
\end{array}
$$
Note that the two sections $X_0 \otimes_R k$ and $X_1\otimes_R k$ 
coincide at the special fibre $H^0(\EE_0, 2\Theta_{|\EE_0})
\cong H^0(\EE,2\Theta) \otimes_R k$.

\section{Equations of $\tilde{V}$ for non-ordinary $X$}

\subsection{Specializing an ordinary curve}

Let $X/k$ be a smooth genus $2$ curve with Hasse-Witt invariant 
equal to $1$. Following e.g. \cite{lange} $X$ is birational
to an affine curve given by an  equation of the form
$$ y^2 + xy = \lambda x^5 + \mu x^3 + x, $$
with $\lambda, \mu \in k$ and $\lambda \not= 0$. The projection
$(x,y) \mapsto x$ defines a separable double cover
$X \ra \pp^1_k$ ramified at $0$ and $\infty$. Let $\pp^1_R$
be the projective line over $R=k[[s]]$ with affine coordinate $x$.
We introduce the family $\XX \ra \pp^1_R$ defined by the projective 
closure of the affine curve with equation
$$ y^2 + (sx^2 +x) y = \lambda x^5 + \mu x^3 +x. $$
The special fibre $\XX_0/k$ equals $X$ and 
the generic fibre $\XX_\eta /K$  of the family $\XX$ is 
a smooth ordinary curve of genus $2$, which is 
birational to the curve (defined over a finite extension of $K$)
given by the standard equation \eqref{affmodbhosle} with
coefficients
\begin{equation} \label{coeffs}
 a= \lambda/s^3, \  b = \alpha^2 + \alpha, \  c =s, \qquad \text{and}
\qquad \alpha^2 = \lambda/s^3 + \mu/s + s.
\end{equation}

Let $\JJ\XX$ be the associated Jacobian scheme and $\JJ\XX[2]/R$ be the
group scheme of $2$-torsion points. Then we have the following
isomorphisms
$$ \JJ\XX[2]_\eta \cong \left(\ZZ/2\ZZ \right)^2 \times \mu_2^2/K, \qquad
\JJ\XX[2]_0 \cong JX[2] \cong \left(\ZZ/2\ZZ \right) \times \mu_2 
\times G_{l,l}/k,$$
where $G_{l,l}$ is the unique self-dual local-local group scheme of 
dimension 1 and length 4. Note that $G_{l,l}$ is isomorphic
to the group scheme of $2$-torsion points $E^{ss}[2]$ (see
section 4.2.2). The \'etale parts of both fibres can be described 
in terms of Weierstrass points as follows.

\bigskip

The $3$ Weierstrass points of $\XX_\eta \ra \pp^1_K$ are $0_\eta$, 
$\infty_\eta$, and
$1_\eta$ with affine coordinate $1/s$, which specialize to $0$, $\infty$
and $\infty$ respectively. We obtain by \eqref{twotorsion} the three
nonzero elements of $\JJ\XX[2]_\eta^{et}$, which we denote by
$[0]_\eta$, $[1]_\eta$ and $[\infty]_\eta$. At the special fibre 
the nonzero 2-torsion point in $JX[2]^{et} \cong \ZZ/2\ZZ$
equals $AJ(0+ \infty)$, which we denote by $[1]_0$. We see that 
$[1]_\eta$ and $[\infty]_\eta$ specialize to $[1]_0 \in 
JX[2]^{et}$, and $[0]_\eta$ specializes to $0$.

\subsection{Decomposing Heisenberg groups}

We are interested in the linear action of the Heisenberg group scheme $\HH/R$,
which is a central extension (see \cite{mum} page 221)
\begin{equation} \label{heis}
0 \lra \mu_2 \lra \HH \lra \JJ\XX[2] \lra 0,
\end{equation}
on the free $R$-module $\WW = H^0(\JJ\XX[2], 2\Theta)$ of rank $4$. We 
choose the splitting over $R$ of the connected-\'etale exact
sequence
$$ 0 \lra \JJ\XX[2]^0 \lra \JJ\XX[2] \lra \JJ\XX[2]^{et} \lra 0 $$
determined by the nonzero $2$-torsion point 
$[1] =  AJ(0+\infty) \in \JJ\XX[2]$.
Note that $\JJ\XX[2]^{et} \cong \ZZ/2\ZZ$ and that $[\infty] \in \JJ\XX[2]$ 
determines a different splitting. Passing to the 
Cartier dual we obtain a decomposition over $R$,
$$ \JJ\XX[2] = \ZZ/2\ZZ \times \mu_2 \times \JJ\XX[2]^{00}.$$
Pulling-back the central extension \eqref{heis} by the
canonical inclusions of $\ZZ/2\ZZ \times \mu_2$ and $\JJ\XX[2]^{00}$
into $\JJ\XX[2]$ we obtain the two Heisenberg groups
$\HH^{et}$ and $\HH^0$
$$ 0 \lra \mu_2 \lra \HH^{et} \lra \ZZ/2\ZZ \times \mu_2 \lra 0,
\qquad
0 \lra \mu_2 \lra \HH^0 \lra \JJ\XX[2]^{00} \lra 0.$$
It is clear that the Heisenberg group scheme $\HH$ \eqref{heis} is
isomorphic to the quotient $\HH^0 \times \HH^{et}/\mu_2$, where
$\mu_2$ acts diagonally on $\HH^0 \times \HH^{et}$. Let 
$\WW^{et}$ and $\WW^0$ be the sub-$R$-modules of $\WW$ fixed by the
subgroups $\HH^0$ and $\HH^{et}$ of $\HH$. By general
theory of Heisenberg groups, $\WW$ is the unique irreducible
$\HH$-module of weight $1$, which implies an $\HH$-isomorphism
$$ \WW \cong \WW^0 \otimes \WW^{et}.$$
Moreover $\WW^0$  (resp. $\WW^{et}$) is the unique  irreducible $\HH^0$
(resp. $\HH^{et}$)-module of weight $1$.

\bigskip

Let $H/k$ be the Heisenberg group scheme associated to 
$\ZZ/2\ZZ \times \mu_2$ ($\cong E[2]$ for any ordinary elliptic curve
$E/k$) and let $W$ be the unique irreducible $H$-module of weight $1$. 
Note that  $W$ is isomorphic (as $H$-module) to the space \eqref{thetaell}. 
It is clear that we have the following isomorphisms  
$$ \HH^{et} \cong H \otimes_k R , \qquad \WW^{et} \cong  W \otimes_k R .$$
We will denote by $\{Z_0,Z_1\}$ the ``constant'' $R$-basis of $\WW^{et}$
induced by the Theta basis $\{ X_0, X_1 \}$ of $W$,i.e.,
$Z_i := X_i \otimes_k 1$.

\bigskip

The structure of the $\HH^0$-module $\WW^0$ is determined by 
analyzing the group scheme $\JJ\XX[2]^{00}$. In section 4.1 we
considered
the $2$-divisible group $\JJ\XX(2)^{00}$ and we showed 
the existence of an elliptic curve $\EE_\XX/R$ such that
$\EE_\XX(2) \cong \JJ\XX(2)^{00}$. In particular $\EE_\XX[2] 
\cong \JJ\XX[2]^{00}$. We observe that the $j$-invariants 
of the  elliptic curves $\EE_\XX/R$ and $\EE/R$ (section
4.2.2) lie in the maximal ideal of $R$, because 
$(\EE_\XX)_0 \cong \EE_0 \cong  E^{ss}$. Therefore there
exists a relation of the form
\begin{equation} \label{unifpar}
s^n = u t^m
\end{equation}
between the two uniformizing parameters
$s$ in \eqref{coeffs} and $t$ in \eqref{eqcubfam}, with $u$ invertible
in $k[[t]]$ and $n,m \in \NN^*$, and we can assume, after passing
to the ramified cover given by \eqref{unifpar}, that the $\HH^0$-module
$\WW^0$ equals $H^0(\EE,2\Theta)$.

\bigskip

In order to have a consistent notation we denote the $R$-basis $\{U,Z\}$ of 
$\WW^0 = H^0(\EE,2\Theta)$ by
$\{Z_0,Z_1\}$ and recall from section 4.2.2 the transition
formulae 
\begin{equation}
X_0  =  t^4 Z_0 + Z_1, \qquad
X_1  =  (1+t^3 + t^4 + t^5)Z_1.  
\end{equation} 
Let $\{x_i\}$ and $\{ z_i \}$ denote the dual $K$-bases of
$\{X_i \}$ and $\{ Z_i \}$ in both spaces $\WW^0$ and
$\WW^{et}$. Then the $4$ tensors $z_{ij} : = z_i \otimes z_j
\in \WW^*$ form an $R$-basis and the dual Theta functions
$x_{ij} := x_i \otimes x_j$ can be expressed as follows
(after chasing denominators)
\begin{equation} \label{trans}
\begin{array}{cc}
x_{00} = (1+t^3 + t^4 + t^5) z_{00}, \qquad  & x_{10} = z_{00} + t^4 z_{10}, \\
x_{01} = (1+t^3 + t^4 + t^5) z_{01}, \qquad  & x_{11} = z_{01} + t^4 z_{11}. 
\end{array}
\end{equation}
Note that the 
coordinate $x_{10}$ specializes to $x_{00}$ and $x_{11}$ specializes
to $x_{01}$. Via the level 2 structure \eqref{leveltwopsi} this parallels
the specialization of the $2$-torsion points $[0]_\eta$, $[1]_\eta$,
and $[\infty]_\eta$ (section 5.1).

\subsection{Specializing the quadrics}

It can be shown as in \cite{lp} section 5 that the identification
$\MM_X \lra |2\Theta|$ extends to the relative case $\XX \lra 
\Spec(R)$,i.e., we have an isomorphism $\MM_\XX \lra \pp(\WW)$
over $\Spec(R)$. Therefore the relative Frobenius morphism 
$\XX \lra \XX_1$ (over $\Spec(R)$) induces by pull-back a rational map 
$$
\begin{array}{lcccr}
\pp(\WW_1) & & \map{\tilde{\VV}} & & \pp(\WW) \\
 & \searrow  & & \swarrow & \\
 &   & \Spec(R) & & 
\end{array}
$$
with $\WW_1 = H^0(\JJ\XX_1, 2\Theta_1)$. We recall that
the map $\tilde{\VV}$ is given by a linear system of $4$ quadrics.
Over the generic point $\eta \in \Spec(R)$ the equations of the map
$\tilde{\VV}_\eta : \pp(\WW_1)_\eta \lra \pp(\WW)_\eta$ are of the
form (\cite{lp} Proposition 3.1 (3))
\begin{equation} \label{eqtildev}
\tilde{\VV} \ : \ x = (x_{ij}) \lms (\lambda_{00} P_{00}(x):
\lambda_{01} P_{01}(x): \lambda_{10} P_{10}(x): \lambda_{11} P_{11}(x))
\end{equation}
with
$$
\begin{array}{l}
P_{00}(x) = x^2_{00} + x^2_{01} + x^2_{10} + x^2_{11}, \\ 
P_{01}(x) = x_{00}x_{01} + x_{10}x_{11},
\end{array}
\qquad
\begin{array}{l}
P_{10}(x) = x_{00}x_{10} + x_{01}x_{11}, \\
P_{11}(x) = x_{00}x_{11} + x_{10}x_{01}.
\end{array}
$$
Here we use the $K$-basis of Theta coordinates $x_{ij}$ on 
$\pp(\WW_1)_\eta$ and $\pp(\WW)_\eta$. Moreover 
Proposition 3.1 relates the coefficients $\lambda_{ij}$
(defined up to a scalar) to the coefficients $a,b,c \in K$ 
\eqref{coeffs},
$$ (\lambda_{00}: \lambda_{01} : \lambda_{10} : \lambda_{11} ) =
(\sqrt{abc}: \sqrt{b}: \sqrt{c}: \sqrt{a}). $$
Since the Theta coordinates $x_{ij}$ are no longer independent after
specialization, we express the equations \eqref{eqtildev}
of $\tilde{\VV}$ in the $R$-basis $\{ z_{ij} \}$ using the transition
formulae \eqref{trans}. A straight-forward computations shows that the
map
\begin{equation} \label{mapr}
\tilde{\VV} \ : \ z = (z_{ij}) \lms (R_{00}(z):
R_{01}(z): R_{10}(z): R_{11}(z))
\end{equation}
is given by the quadrics
\begin{equation} \label{quadricsrij}
\begin{array}{rcl}
R_{00}(z) &  = & \frac{\sqrt{abc}}{1+t^3 + t^4 + t^5} 
\left[ (t^{12} + t^{16} + t^{20})
(z^2_{00} + z^2_{01}) + t^{16} (z^2_{10} + z^2_{11})  \right], \\
R_{01}(z) &  = & \frac{\sqrt{b}}{1+t^3 + t^4 + t^5} \left[ 
(t^{12} + t^{16} + t^{20})
z_{00}z_{01} + t^8(z_{00}z_{11} + z_{10}z_{01}) + t^{16} z_{10}z_{11}
\right], \\
R_{10}(z) & = & \frac{1}{t^4} \left[ R_{00} + \sqrt{c} (1+t^6 + t^8 + t^{10})
\left( z^2_{00} + z^2_{01}  + t^8 (z_{00}z_{10} + z_{11} z_{01}) \right)
 \right], \\
R_{11}(z) & = & \frac{1}{t^4} \left[ R_{01} + \sqrt{a} 
(1+t^6 + t^8 + t^{10})t^8 (z_{00}z_{11} + z_{01}z_{10}) \right].
\end{array}
\end{equation}
Note that we square the coefficients in \eqref{trans} when
considering coordinates on $\pp(\WW_1)$. At the special fibre
the map \eqref{mapr} specializes to $\tilde{\VV}_0$ obtained
by putting $t=0$ after having divided the $R_{ij}$'s by $t^\alpha$,
where $\alpha$ is the lowest valuation appearing in the
expressions \eqref{quadricsrij}. The map $\tilde{\VV}_0$
coincides with the Verschiebung $\tilde{V}$ of the curve $X$,
because both maps extend to rational maps over $R$ and
coincide over $K$. Since the image of $\tilde{V}$ is non-degenerate
(it contains the Kummer surface $\kum_X \subset |2\Theta|$), the
lowest valuations for each of the $4$ quadrics $R_{ij}$ coincide
(otherwise the image of $\tilde{\VV}_0$ is contained in a
hyperplane).

\bigskip

We work out the specialization of the quadrics as follows. We write
$\nu = \frac{m}{n}$ and replace $s$ by $v t^\nu$, with $v \in
k[[t]]$ invertible, in the expression of the coefficients $a,b,c$
of \eqref{coeffs}. Note that the (rational) valuations of
$\sqrt{a},\sqrt{b},\sqrt{c}$ are $-\frac{3}{2}\nu,-\frac{3}{2}\nu,
\frac{1}{2}\nu$ respectively. First we observe
that the valuations of $R_{01}$ and $R_{00}$ equal $8-\frac{3}{2}\nu$ and
$12-\frac{5}{2}\nu$ respectively. Since they coincide, 
we obtain $\nu = 4$,i.e.,
$$ R_{00} = t^2 (z^2_{00} + z^2_{01}) + \ \text{h.o.t.}, \qquad
R_{01} = t^2 (z_{00}z_{11} + z_{10}z_{01}) + \ \text{h.o.t.}, $$
up to some multiplicative nonzero constants. Next we see that the 
expansions of $R_{10}$ and $R_{11}$ are given by
$$ R_{10} = t^2 (z^2_{00} + z^2_{01} + z^2_{10} + z^2_{11}) + \ \text{h.o.t.},
\qquad
R_{11} = t^2 (z_{00}z_{01}) + \ \text{h.o.t.}, $$
up to some multiplicative nonzero constants and some multiple of
$R_{00}$ and $R_{01}$ respectively. Thus we have shown

\begin{thm} \label{mainthm}
Let $X$ be a smooth curve with Hasse-Witt invariant equal to $1$. 
There exist coordinates $\{ z_{ij} \}$ on 
$|2\Theta_1|$ and $\{ y_{ij} \}$ on $|2\Theta|$ such that the equations 
of $\tilde{V}$ are given
by 
$$ |2\Theta_1| \map{\tilde{V}} |2\Theta|, \qquad
z=(z_{ij}) \lms y = (y_{ij}) =  
(\lambda_{00} Q_{00}(z) : \lambda_{01} Q_{01}(z) :
\lambda_{10} Q_{10}(z) : \lambda_{11} Q_{11}(z) ) $$
with 
$$ 
\begin{array}{l}
Q_{00}(z) = z^2_{00} + z^2_{01}, \\ 
Q_{01}(z) = z_{00}z_{11} + z_{10}z_{01},
\end{array}
\qquad
\begin{array}{l}
Q_{10}(z) = z^2_{00} + z^2_{01} + z^2_{10} + z^2_{11}, \\
Q_{11}(z) = z_{00}z_{01},
\end{array}
$$
and the $\lambda_{ij}$'s are nonzero constants depending on the
curve $X$.
\end{thm}

\begin{rem}
We note that the equations of $\tilde{V}$ given in Theorem 
\ref{mainthm} are written in two sets of coordinates on $|2\Theta|$
and $|2\Theta_1|$ which do not necessarily correspond under
the $k$-semi-linear isomorphism $JX_1 \ra JX$.
\end{rem}

\begin{rem}
In case $X$ is a non-ordinary curve with Hasse-Witt invariant 
equal to $0$, i.e., $X$ is supersingular, we observe that the $2$-divisible
group $JX(2)= JX(2)^{00}$ (see section 4.1) is self-dual, of dimension $2$
and height $4$. There exists a finite number of isomorphism classes
of such $2$-divisible groups over $k$ (see \cite{de} page 93). Moreover 
one can show that $JX(2)$ cannot be isomorphic to the product
$E^{ss}(2) \times E^{ss}(2)$. 
\end{rem}

\subsection{Applications}

As in \cite{lp} section 6 we can easily deduce from Theorem 
\ref{mainthm} a full description of the Verschiebung $\tilde{V}$.
Since the computations are analoguous to those of \cite{lp} 
Proposition 6.1, we leave them to the reader.

\begin{prop} \label{descr}
Let $X$ be a smooth genus $2$ curve with Hasse-Witt invariant 
equal to $1$.  
\begin{enumerate}
\item There exists a unique stable bundle $E_{BAD} \in \MM_{X_1}$, which
is destabilized by the Frobenius map, i.e., $F^* E_{BAD}$ is
not semi-stable. We have $E_{BAD} = {F}_* B^{-1}$ and its projective
coordinates are $(0:0:1:1)$.

\item Let $H_1$ be the hyperplane in $|2\Theta_1|$ defined by 
$z_{00} + z_{01} = 0$. The map $\tilde{V}$ contracts $H_1$ to the conic
$\kum_X \cap H$, where $H$ is the hyperplane in $|2\Theta|$ defined by 
$y_{00} = 0$.

\item The fiber of $\tilde{V}$ over a point $[E] \in \MM_X$ is

\begin{itemize}
\item a non-degenerate $\ZZ/2\ZZ$-orbit of a point $[E_1] \in
\MM_{X_1}$, if $[E] \notin H$
\item empty, if $[E] \in H \setminus (H \cap \kum_X)$
\item a projective line passing through $E_{BAD}$, if $[E] \in
H \cap \kum_X$
\end{itemize}

In particular, $\tilde{V}$ is dominant and non-surjective. 
The separable degree of $\tilde{V}$ is $2$.

\end{enumerate}

\end{prop}

\section{Equations of $\tilde{V}$ in characteristic $3$}

Let $X$ be a smooth curve of genus $2$ defined over a field of
characteristic $p = 3$. The main result of this section is

\begin{thm} \label{theochar3}
\begin{enumerate}
\item There exists an embedding $\alpha: \kum_X \hookrightarrow
|2\Theta_1|$ such that the equality of hypersurfaces in $|2\Theta_1|$
$$ \tilde{V}^{-1} (\kum_X) = \kum_{X_1} \cup \alpha(\kum_X) $$
holds {\em set-theoretically}.
\item The cubic equations of $\tilde{V}$ are given by the $4$ partials
of the quartic equation of the Kummer surface $\alpha(\kum_X) \subset
|2\Theta_1|$. In other words, $\tilde{V}$ is the polar map of the
surface $\alpha(\kum_X)$.
\end{enumerate}
\end{thm}

\begin{proof}
Given $L \in JX$ and a section
$\varphi \in H^0(X,\omega \otimes L^2)$, we can consider, using 
adjunction and relative duality for the map $F$, the vector bundle map 
$$ F_* L \map{\varphi} L \otimes \omega, $$
where we also write $L$ for the pull-back $\iota^* L$ under the 
$k$-semi-linear isomorphism $\iota: X_1 \ra X$. We define 
$E_L := \ker(\varphi)$. Hence there is an exact sequence
of vector bundles over $X_1$
\begin{equation} \label{esel}
 0 \lra E_L \lra F_*L \map{\varphi} L \otimes \omega \lra 0. 
\end{equation}
It is staight-forward to show that $E_L$ has determinant $\cO_{X_1}$, 
is semi-stable, that $E_L = E_{L^{-1}}$, (if $L^{-1} \not= L$) 
and that $F^*E_L$ is  $S$-equivalent to $L \oplus L^{-1}$. Indeed,
by adjunction, the inclusion $E_L \ra F_* L$ induces a nonzero
map $F^* E_L \ra L$.

\bigskip

We observe that for $L$ such that $L^2 \not= \cO$, $\dim H^0(\omega L^2)
=1$, hence $\varphi$ is uniquely defined (up to a scalar). For $L$
such that $L^2 = \cO$, we obtain a projective line $\pp^1_L$ of rank $2$ vector
bundles $E_{L,\varphi}$ with $\varphi \in |\omega|$. The variety
of pairs $(L,\varphi)$ is isomorphic to the blowing-up
$\mathrm{Bl}_2(JX)$ of $JX$ at the $16$ $2$-torsion points $L \in
JX[2]$. Hence we obtain a morphism
\begin{equation*}
e: \mathrm{Bl}_2(JX) \lra \MM_{X_1} \cong |2\Theta_1|, \qquad
(L,\varphi) \lms E_L.
\end{equation*} 

We will determine the image of the morphism $e$. First we consider the
image of the $16$ exceptional divisors $\pp^1_L$ for $L \in JX[2]$. We
identify $H^0(X_1,\omega_1) = H^0(X, \omega)$.

\begin{lem}
Given $L \in JX[2]$ and $w$ a Weierstrass point,i.e., $2w \in |\omega|$.
Then we have $E_{L,\varphi} \cong B \otimes L^{-1}(-w)$ with 
$Div(\varphi) = 2w$.
\end{lem}

\begin{proof}
We tensorize the exact sequence \eqref{esel} with $L$, using 
$F_*L \otimes L = F_*(L \otimes F^*L) = F_* \cO _X$ (projection
formula)
$$
\begin{array}{ccccccccc}
 & & & & \cO_{X_1} & & & & \\
 & & & & \downarrow & \searrow^\varphi & & & \\
0 & \lra & E_{L,\varphi} \otimes L & \lra  & F_* \cO_X & \lra &  \omega_{X_1} &
\lra & 0 \\
 & & & \searrow^{\hat{\varphi}} & \downarrow & & & & \\
 & & & & B & & & & 
\end{array}
$$
The vertical arrows form an exact sequence (see \cite{ray} section 4.1).
The upper diagonal map defined as composite map $\cO_{X_1} \ra 
F_* \cO_X \ra \omega_{X_1}$ equals $\varphi$,
which implies that the lower diagonal map $\hat{\varphi} : 
E_{L,\varphi} \otimes L \ra B$ vanishes at $w$. Hence, using stability of 
$B$ and $\det B = \omega$, we obtain $E_{L,\varphi} \otimes L \cong B(-w)$.
\end{proof}

We observe that $L(w)$ is a theta-characteristic $\kappa$. Moreover the
image of $\pp^1_L$ in $|2\Theta_1| = \pp^3$ is a conic and it is
straight-forward to establish that the $16$ points $B \otimes \kappa^{-1}$
(with $\kappa$ varying over the set of theta-characteristics) and
the $16$ hyperplanes spanned by the $16$ conics $e(\pp^1_L)$ form
a classical $16_6$-configuration of nodes and tropes, whose
associated Kummer surface equals $\kum_X$ (see \cite{gd}). We 
recover the curve $X$
e.g. by taking the $2:1$ cover over a $\pp^1_L$ branched at the $6$
nodes on $\pp^1_L$. Since $\deg e(\mathrm{Bl}_2(JX)) =4$ and 
$e(\mathrm{Bl}_2(JX))$ is singular at the $16$ points $B \otimes
\kappa^{-1}$, the image $e(\mathrm{Bl}_2(JX))$ coincides with 
$\alpha(\kum_X)$ for an embedding $\alpha: \kum_X \hookrightarrow
\MM_{X_1}$

\bigskip

By Proposition \ref{degreeV} the map $\tilde{V}$ is given by a linear system
$|\LL|$ of $4$ cubics on $|2\Theta_1|$. The key fact underlying
Theorem \ref{theochar3} is a striking relationship between cubics
and quartics on $|2\Theta_1|$ (\cite{vg} Proposition 2): the $4$
cubics in $|\LL|$ are the $4$ partial derivatives of a Heisenberg
invariant quartic, whose zero divisor we denote by $Q \subset
|2\Theta_1|$.

\bigskip

By \cite{ray} Remarque 4.1.2 (2) the $16$ points $B \otimes \kappa^{-1}$
(a $JX_1[2]$-orbit) are contained in the base locus of $|\LL|$. Hence
$Q$ is singular at the $16$ points $B \otimes \kappa^{-1}$. Since a
Kummer surface is determined by its singular points, we conclude 
that $Q \cong \alpha(\kum_X)$ and that $|\LL|$ is given by the $4$
partials of $\alpha(\kum_X)$. This completes the proof of Theorem
\ref{theochar3}.
\end{proof}

\begin{rem}
\begin{enumerate}
\item The rational map $e: \kum_X \lra \MM_{X_1}$ (defined away from
the $16$ nodes of $\kum_X$) is the birational inverse of $\tilde{V}$.

\item One has the following scheme-theoretical equality (among divisors
in $|2\Theta_1|$)
$$ \tilde{V}^{-1}(\kum_X) = \kum_{X_1} + 2\alpha(\kum_X).$$
\end{enumerate}
\end{rem}

\begin{cor} \label{corchar3}
\begin{enumerate}
\item The map $\tilde{V}$ has exactly $16$ base points, which
correspond bijectively to the
\begin{itemize}
\item $16$ nodes of the surface $\alpha(\kum_X) \subset |2\Theta_1|$
\item $16$ stable rank $2$ vector bundles $B \otimes \kappa^{-1} \in
\MM_{X_1} \cong |2\Theta_1|$, where $B$ is the bundle of locally
exact differentials and $\kappa$ a theta-characteristic of $X_1$. 
\end{itemize}

\item The map $\tilde{V}$ is surjective, separable and of degree $11$.
\end{enumerate}
\end{cor}

\begin{proof}
We only have to show part 2, since part 1 is clear from the proof of
Theorem \ref{theochar3}. We recall that the rational map 
$\tilde{V}$, which is
defined away from the $16$ points $B \otimes \kappa^{-1}$, coincides
with the polar map of the Kummer surface $\kum_X$. It is well-known
that $\tilde{V}$ can be resolved into a morphism $\mathcal{V}:
\mathrm{Bl}(|2\Theta_1|) \lra |2\Theta|$ by blowing-up these
$16$ points in $|2\Theta_1|$. We denote by $E_\kappa$ the 
exceptional divisor over $B \otimes \kappa^{-1}$ and by $H_\kappa
\subset |2\Theta|$ the hyperplane $\mathcal{V}(E_\kappa)$. Note that
the $H_\kappa$ are the $16$ tropes of $\kum_X \subset |2\Theta|$ and that
$\mathcal{V}_{|E_\kappa}$ is a linear isomorphism. It is clear
that the image of $\tilde{V}$ contains the complement of the $16$
hyperplanes $H_\kappa$.

\bigskip

Let us check that the $H_\kappa$ are also contained in the image of
$\tilde{V}$: a simple computation shows that the cubic $C_\kappa :=
\tilde{V}^{-1}(H_\kappa) \subset |2\Theta_1|$ is singular at the point
$B \otimes \kappa^{-1}$ and that the restriction of $\tilde{V}$ to the
cubic $C_\kappa$ concides with the (birational) projection with center
$B \otimes \kappa^{-1}$. Moreover the projectivized tangent cone
at $B \otimes \kappa^{-1}$  to $C_\kappa$ is the conic $Q_\kappa
\subset H_\kappa$ through the $6$ nodes (recall that $2Q_\kappa =
H_\kappa \cap \kum_X$). Therefore any point in $H_\kappa \setminus
Q_\kappa$ lies in the image of $\tilde{V}$. To finish the argument
we observe that $Q_\kappa \subset \kum_X$ and that $\tilde{V}:
\kum_{X_1} \lra \kum_X$ is surjective.
\end{proof}

\begin{rem}
\begin{enumerate}
\item We recall (\cite{lp} Remark 6.2) that surjectivity only holds for 
$S$-equivalence classes (not isomorphism classes!). In fact, there 
always exist semi-stable bundles $E$ which do not descend 
by Frobenius.

\item The number of base points and the degree of $\tilde{V}$ was 
also obtained in \cite{o} by computing the number of connections
(on certain unstable bundles) with zero $p$-curvature. 

\item It would be interesting to have an explicit description
of the $11$ vector bundles in a general fiber $\tilde{V}^{-1}(E)$ of the
polar map $\tilde{V}$.

\end{enumerate}
\end{rem}
 
\section{Appendix: base points of $\tilde{V}$}  

In this section we consider a smooth curve $X$ of genus $g \geq 2$
defined over an algebraically closed field $k$ of caracteristic $p>0$.
We denote by $\MM_X(r)$ (resp. $\MM_{X_1}(r)$) the moduli
space of semi-stable rank $r$ vector bundles over $X$ (resp. $X_1$)
and by $\LL$ (resp. $\LL_1$) the determinant line bundle over 
$\MM_X(r)$ (resp. $\MM_{X_1}(r)$). The relative Frobenius
map $F: X \ra X_1$ induces by pull-back a rational map
$$ V: \MM_{X_1}(r) \lra \MM_X(r), $$
called the Verschiebung. Let $\II$ be the indeterminacy locus of $V$, i.e.,
the closed subscheme of $\MM_{X_1}(r)$ consisting of semi-stable
bundles $E$ such that $F^* E$ is not semi-stable.  Let $U = 
\MM_{X_1}(r) \setminus \II$ be the open subset where $V$ is a 
morphism.

\subsection{General facts}

\begin{prop} 
We have an isomorphism $V^*(\LL) \cong \left(\LL^{\otimes p}_{1}
\right)_{|U}$.
\end{prop}

\begin{proof}
Let $\M_X(r)$ and $\M_{X_1}(r)$ be the moduli stacks of
rank $r$ vector bundles over $X$ and $X_1$ and let $\EE$ and
$\EE_1$ be the universal bundles with trivialized determinant
on $X \times \M_X(r)$ and $X_1 \times \M_{X_1}(r)$. It is
well-known that the inverses of the determinant of 
cohomology, which we denote by $\det Rp_* \EE$ and $\det Rp_{1*}
\EE_1$ descend (after restriction to the semi-stable loci) to
the line bundles $\LL$ and $\LL_1$ on the moduli spaces
$\MM_X(r)$ and $\MM_{X_1}(r)$. Now, since $\det Rp_*$ 
commutes with base change, we have an isomorphism
over the moduli stack $\M_{X_1}(r)$
$$ V^*(\det Rp_* \EE) \cong \det Rp_* \left( 
(F\times \mathrm{Id})^* \EE_1 \right).$$
Moreover we have a commutative diagram
\begin{eqnarray*}
X \times \M_{X_1}(r) & \map{F \times \mathrm{Id}} & X_1 \times
\M_{X_1}(r) \\
 \searrow^p &  & \swarrow^{p_1} \\
 & \M_{X_1}(r) &
\end{eqnarray*} 
where $p$ and $p_1$ denote the projections on  the second
factor. Since $F \times \mathrm{Id}$ is an affine morphism,
we have $R^1 (F \times \mathrm{Id})_* = 0$. Hence
$$\det Rp_* \left( (F\times \mathrm{Id})^* \EE_1 \right)  \cong
\det Rp_{1*} \left( (F \times \mathrm{Id})_* (F\times \mathrm{Id})^* 
\EE_1 \right) \cong \det Rp_{1*} (\EE_1 \boxtimes F_* \cO_{X}).$$
The last equality follows from the projection formula. Using 
a filtration by line bundles of the rank $p$ bundle
$F_* \cO_X$ and by showing that $\det Rp_{1*} (\EE_1 \boxtimes 
\cO_{X_1}(D)) = \det Rp_{1*} (\EE_1)$ for an effective divisor $D$,
we show that $\det Rp_{1*} (\EE_1 \boxtimes F_* \cO_X) \cong
(\det Rp_{1*} \EE_1)^{\otimes p}$. We obtain the isomorphism
of the lemma by descent on $U$.
\end{proof}

\begin{prop} \label{degreeV}
If $g=2$ and $r=2$, then $\dim \II = 0$ and the rational map $\tilde{V}$ 
is given by polynomials of degree $p$.
\end{prop}

\begin{proof}
The fact that $\dim \II = 0$ is proved in Theorem 3.2 \cite{jx}. This
implies that $V^*(\LL)$ extends uniquely to $\LL_1^{\otimes p}$ over
$\MM_{X_1}$ and the lemma follows from the isomorphism $\LL_1 \cong
\cO_{\pp^3}(1)$.
\end{proof}

\begin{rem}
For general $g,r,p$ we do not know an estimate of the dimension of
$\II$.
\end{rem}

\subsection{Existence of base points}

\begin{thm} \label{exisbp}
The indeterminacy locus $\II$ is non-empty.
\end{thm}

\begin{proof}
First it will be enough to show non-emptiness of $\II$ in the case 
$r=2$, since taking direct sums with the trivial bundle implies
non-emptiness for arbitrary $r$.  Secondly it suffices to show
non-emptiness of $\II$ after a field extension $k'/k$, with $k'$
algebraically closed. 
\bigskip

Let $\overline{\M_g}$ be the coarse moduli space of stable
genus $g$ curves defined over $k$, which is an irreducible
projective variety \cite{dm}. Let $\eta$ be the generic point
of $\overline{\M_g}$. The choice of a geometric point $\overline{\eta}$
over $\eta$ defines a smooth curve $\XX_{\overline{\eta}}$  over
$\overline{k(\eta)}$, the algebraic closure of the function field
$k(\eta)$ of $\overline{\M_g}$. The curve $\XX_{\overline{\eta}}$ is
defined over a finite extension $K$ of $k(\eta)$ and we denote by 
$\XX_K$ some model of $\XX_{\overline{\eta}}$, i.e., $\XX_K \times_K
\overline{k(\eta)} \cong \XX_{\overline{\eta}}$.

\bigskip

The curve $X/k$ defines a $k$-rational point $x$ of $\overline{\M_g}$,
which lies in the closure of $\eta$. The local ring $A_x$ at the
generic point of the exceptional divisor of the blowing-up
of $\overline{\M_g}$ at the point $x$ is a discrete valuation ring
with fraction field $k(\eta)$ and residue field containing $k$.
By the stable reduction theorem (Corollary 2.7 \cite{dm}) there
exists a finite extension $L$ of $K$, and therefore also of $k(\eta)$,
such that $\XX_L$ is the generic fibre of a stable curve $\XX$
over the integral closure $A$ of $A_x$ in $L$. Note that $A$ is
a discrete valuation ring with fraction field $L$ and with residue
field, denoted by $k(s)$, containing $k$. Moreover the diagram
\begin{eqnarray*}
\Spec(A) & & \\
\downarrow \ \ \ & \searrow^{\XX} & \\
\Spec(A_x) & \lra & \overline{\M_g} 
\end{eqnarray*}
commutes when restricted to $\Spec(L) \hookrightarrow \Spec(A)$ and
therefore commutes because $\overline{\M_g}$ is separated. It follows
that the special point $s \in \Spec(A)$ maps to $x$, i.e., there
exists an isomorphism $X \times_k \overline{k(s)} \cong 
\XX_s \times_{k(s)} \overline{k(s)}$.

\bigskip

In summary, we have constructed a stable curve $\XX$ over a discrete
valuation ring $A$  with generic fibre $\XX_L$ and geometric
special fibre isomorphic to $X \times_k \overline{k(s)}$. The 
fraction field of $A$ is $L$ and its residue field $k(s)$.

\bigskip

We now choose a tree of $\pp^1_k$'s, denoted by $X'$, 
defining a closed point $x'$
in the boundary of $\overline{\M_g}$. Repeating the above contruction
with $x'$ instead of $x$, we obtain a stable curve $\XX'$ over a
discrete valuation ring $A'$ with generic fibre $\XX_{L'}$ and
geometric special fibre isomorpic to $X' \times_k \overline{k(s')}$.
The fraction field of $A'$ is $L'$, a finite extension of $k(\eta)$, 
and its residue field is $k(s')$. Moreover the isomorphism
$X' \times_k \overline{k(s')} \cong \XX_{s'} \times_{k(s')} \overline{k(s')}$
is defined over a finite extension of $k(s')$.

\bigskip

We choose a finite extension of $k(\eta)$ containing both
$L$ and $L'$, which we call again $L$, and take the 
integral closures in $L$ of $A$ and
$A'$, which we call again $A$ and $A'$. 
Thus we have constructed two stable curves
$\XX$ and $\XX'$ over $A$ and $A'$ such that $\XX_L \cong \XX'_L$ and which
specialize to $X$ and $X'$ respectively.   

\bigskip
Let $\hat{L}$ be the fraction field of the completion $\hat{A'}$ of $A'$. By
construction the curve $\XX_{\hat{L}} \cong \XX'_{\hat{L}}$ is a
Mumford-Tate curve and, by the main result of \cite{g}, there exists
a stable rank $2$ vector bundle $\hat{\EE}$ over $\XX_{\hat{L}}$ such 
that $F^* \hat{\EE}$ is not semi-stable.

\begin{lem}
There exists a finite extension $L_1$ of $L$ contained in the field
$\hat{L}$ and a stable bundle $\EE_1$ over $\XX_{L_1}$ such that
$\EE_1 \otimes_{L_1} \hat{L} \cong \hat{\EE}$ and $F^* \EE_1$ is not
semi-stable.
\end{lem}   

\begin{proof}
Let $\hat{\pi} : F^* \hat{\EE} \ra \hat{\LL}$ be a maximal 
destabilizing quotient of $\hat{\EE}$. There exist a models
$\EE_{k(S)}, \LL_{k(S)}$ and $\pi_{k(S)}$ of $\hat{\pi}$ over
$\XX_{k(S)}$, where $k(S)$ is an extension of finite type of $L$. The
field $k(S)$ is the function field for some algebraic variety $S$
over $L$. Shrinking $S$ if necessary, one can assume that 
$\pi_{k(S)}$ comes from 
$$\pi_S : F^* \EE_S \ra \LL_S,$$   
where $\EE_S$ is a family of stable bundles over $\XX_L$ parametrized
by $S$ (stability is an open condition). We now choose a closed point $s \in S$
and pull-back the family $\EE_S$ under the inclusion $s \hookrightarrow
S$. We thus obtain a stable bundle $\EE_{L_1}$ over $\XX_{L_1}$, where
$L_1$ is the residue field at the point $s$, which is a finite 
extension of $L$. 
\end{proof}

Again we take the integral closures $A_1$ and $A'_1$ of the
discrete valuation rings $A$ and $A'$ in $L_1$. By the previous
lemma we have a stable bundle $\EE_1$ and a destabilizing quotient
$\LL_1$ over $\XX_{L_1}= \XX'_{L_1}$
$$\pi_{L_1} : F^* \EE_1 \ra \LL_1.$$
After possibly taking a finite extension of $L_1$, we can assume
\cite {langton} that
$\EE_1$ and $\LL_1$ have models over $\XX \ra \Spec(A)$
with $\left(\EE_1\right)_{\overline{s}}$ 
semi-stable over $X \times_k \overline{k(s)}$.
By semi-continuity, we have
$$ \mathrm{Hom} (F^* \EE_{1\overline{s}}, \LL_{\overline{s}}) \not = 0,$$
which shows that  $F^* \EE_{1\overline{s}}$ is not semi-stable.

\end{proof}

\bigskip

\flushleft{Yves Laszlo \\
Universit\'e Pierre et Marie Curie, Case 82 \\
Analyse Alg\'ebrique, UMR 7586 \\
4, place Jussieu \\
75252 Paris Cedex 05 France \\
e-mail: laszlo@math.jussieu.fr}

\bigskip

\flushleft{Christian Pauly \\
Laboratoire J.-A. Dieudonn\'e \\
Universit\'e de Nice Sophia Antipolis \\
Parc Valrose \\
06108 Nice Cedex 02 France \\
e-mail: pauly@math.unice.fr }

\end{document}